\documentclass[12pt,reqno,oneside]{amsart}
 \usepackage{verbatim,amsmath,amssymb,cite,epsfig,xspace}
\usepackage{color,cite,graphicx}

\numberwithin{equation}{section}

\theoremstyle{plain}
\newtheorem{theorem}{Theorem}[section]
\newtheorem{lemma}{Lemma}[section]
\newtheorem{corollary}{Corollary}[section]
\theoremstyle{definition}
\newtheorem{definition}{Definition}[section]

\theoremstyle{remark}
\newtheorem{remark}{Remark}[section]

\begin{document}
\title[Minimal blow-up initial data]{Minimal $L^3$-initial data for potential Navier-Stokes singularities}
\author{Hao Jia, Vladim\'{i}r \v{S}ver\'{a}k}

\begin{abstract}

We give a simple proof of the existence of initial data with minimal $L^3$-norm for potential Navier-Stokes singularities, recently established in ``Gallagher, I., Koch, G.S., Planchon, F., \emph{A profile decomposition approach to the $L^{\infty}_t(L^3_x)$ Navier-Stokes regularity criterion}, arXiv:1012.0145v2'' with techniques based on profile decomposition. Our proof is more elementary, and is based on suitable splittings of initial data and energy methods. The main difficulty in the $L^3$ case is the lack of compactness of the imbedding $L^3_{\rm loc}\hookrightarrow L^2_{\rm loc}$.
\end{abstract}

\maketitle


\begin{section}{Introdution}
We consider the initial value problem for Navier-Stokes equations (NSE) in $R^3$:
\begin{eqnarray}\label{eq:no1}
 \left.\begin{array}{rl}
    \partial_t u-\Delta u+u\cdot\nabla u+\nabla p&=0\\
                                          \mbox{div}~~u&=0
    \end{array}\right\}\quad &&\mbox{in $R^3\times (0,\infty)$,}\\
 u(\cdot,0)=u_0 \quad &&\mbox{in $R^3$.}
\end{eqnarray}
It is well known that if divergence free $u_0$ belongs to one of a number of `critical' spaces with respect to the natural scaling
\begin{eqnarray*} 
u(x,t)&\longrightarrow& \lambda u(\lambda x,\lambda^2 t) {\rm ~~~for~~} \lambda>0,\\
u_0(x)&\longrightarrow& \lambda u_0(\lambda x) {\rm ~~~for~~} \lambda>0,
 \end{eqnarray*}
such as $L^3(R^3)$ and $\dot{H}^{1/2}(R^3)$, NSE has a unique local `mild solution' (see \cite{KoTa, Kato} and references therein). It is not clear whether such mild solutions exist for all time or whether singularities might develop in finite time. In \cite{RuSv} it was shown that there exists a minimal blowup initial data in $\dot{H}^{1/2}(R^3)$ assuming that some initial data in $\dot{H}^{1/2}$ would produce finite time singularity. A natural question is if the result of \cite{RuSv} can be extended to the $L^3(R^3)$ setting. The main tools of \cite{RuSv} are the stability of singularities and compactness for a sequence of suitable weak solutions uniformly bounded in energy norm, certain estimates of the so called `Leray solutions' together with a uniqueness theorem of Leray solutions when we have a `good solution'. One of the crucial points in the proof of uniqueness is the compactness of the embedding $\dot{H}^{1/2}_{loc}(R^3)\hookrightarrow L^2_{loc}(R^3)$. In the case of $L^3(R^3)$ we lose this compactness. Thus the question is whether one can avoid using the compactness. This was done in \cite{IGF}, by using the technique of profile decompositions. Here we present another way to overcome the difficulty.  Moreover, we also recover the compactness of the set of `minimal blow up initial data' in $L^3(R^3)$ modulo translations and scalings, which is not established in \cite{IGF} (see section 5 below).  Our main tool is the following simple observation:
\begin{lemma}
Let $u$ be a Leray solution with divergence free initial data $u_0\in L^3(R^3)$. Then there exists a nonnegative function $h(t)$ depending only on $\|u_0\|_{L^3(R^3)}$, such that $\lim_{t\to 0+}h(t)=0$ and 
\begin{equation}
\|u(\cdot,t)-e^{\Delta t}u_0\|_{L^2(B_1(x_0))}\leq h(t) {\rm ~~for~~ any~~} x_0\in R^3, {\rm ~~a.e.~~} 0\leq t<1.
\end{equation}
\end{lemma}
The proof is based on a suitable splitting of the solutions.  The important point is $h(t)$ in the lemma only depends on the $L^3$-norm of initial data, which gives a certain uniformity of strong continuity in $L^2_{loc}(R^3)$ at time $0$ for a sequence of solutions with initial data uniformly bounded in $L^3(R^3)$. This lemma is, in fact, already sufficient to extend the arguments in \cite{RuSv} to the $L^3$ case.
We present in some detail proofs of certain estimates and uniqueness results for Leray solutions which were proved in Lemari\'{e}-Rieusset \cite{LemRie} and used in\cite{RuSv}, since in the situation considered here the proofs significantly simplify. We will often refer to \cite{RuSv}, as the general ideas are similar and we provide more detailed proof of some points which in \cite{RuSv} were only sketched, and vice versa. We also refer the reader to a recent paper~\cite{GS1} where a related method of using comparisons with solutions of the linear problem is used.

\medskip
\noindent
\textbf{Notation:} We will denote $C$ as an absolute positive number, $C(\alpha,\lambda,\dots)$ denote a positive constant depending on $\alpha,~~\lambda$ and so on. We adopt the convention that nonessential constants $C$ may change from line to line. $B_r(x_0)\subset R^3$ means a ball with radius $r$ centered at $x_0$. $Q_r(x_0,t_0):=B_r(x_0)\times (t_0-r^2,t_0)\subset R^3\times R$, and $Q_r:=Q_r(0,0)$. For vectors $a,~~b$, $a\otimes b$ is a matrix with $(a\otimes b)_{ij}=a_ib_j$. For two matrices $a,~~b$, $(a:b):=a_{ij}b_{ij}$ where we assume the usual Einstein summation convention. We will use $u_0$ as a divergence free initial data for NSE, unless defined otherwise.
\end{section}

\begin{section}{Leray solutions}
In \cite{Leray} J.Leray showed, among many other important results, the existence of a globally defined weak solution $u(x,t)$ to (\ref{eq:no1}) with $u_0\in L^2$ using a priori energy estimates. The regularity and uniqueness of such solutions are open. Later, Calderon \cite{Cald} generalized Leray's theory of weak solutions to the case $u_0\in L^p$. In \cite{LemRie},
Lemari\'{e}-Rieusset constructed global weak solutions with initial data in the space of uniformly locally integrable functions with certain decay at infinity. Here, we recall some results in \cite{LemRie} and present their proofs in some detail for the sake of completeness.
\begin{definition}{(Leray solution)}
$u\in L^2_{loc}(R^3\times [0,\infty))$ is called a Leray solution to NSE with initial data $u_0$ if it satisfies:

\smallskip
\noindent
i)   ${\rm ess}\sup_{0\leq t<R^2}\sup_{x_0\in R^3}\int_{B_R(x_0)}\frac{|u|^2}{2}(x,t)dx+\sup_{x_0\in R^3}\int_0^{R^2}\int_{B_R(x_0)}|\nabla u|^2dxdt<\infty$, and
\begin{equation}\label{decay}
\lim_{|x_0|\to \infty}\int_0^{R^2}\int_{B_R(x_0)}|u|^2(x,t)dxdt=0,
\end{equation}
for any $R<\infty$.

\smallskip
\noindent
ii)  for some distribution $p$ in $R^3\times (0,\infty)$, $(u,p)$ verifies NSE (\ref{eq:no1}) in the sense of distributions and for any compact set $K\subseteq R^3$, $\lim_{t\to 0+}\|u(\cdot,t)-u_0\|_{L^2(K)}=0$.

\smallskip
\noindent
iii)  $u$ is suitable in the sense of Caffarelli-Kohn-Nirenberg, more precisely, the following local energy inequality holds:
\begin{equation}\label{eq:no2}
\int_0^{\infty}\int_{R^3}|\nabla u|^2\phi(x,t)dxdt\leq \int_0^{\infty}\int_{R^3}\frac{|u|^2}{2}(\partial_t\phi+\Delta \phi)+\frac{|u|^2}{2}u\cdot \nabla \phi+pu\cdot\nabla\phi dxdt,
 \end{equation}
for any smooth $\phi\ge 0$ with ${\rm supp}~~ \phi\Subset R^3\times (0,\infty)$. The set of all Leray solutions starting from $u_0$ will be denoted as $\mathcal{N}(u_0)$.
\end{definition}

\noindent
\textbf{Remarks:} In the case the initial data is in $L^2(R^3)$, the notion of Leray-Hopf weak solution is often used (see \cite{LaSe} for example). The difference is that Leray-Hopf weak solution is in $L^{\infty}_tL^2_x\cap L^2_t\dot{H}^{1}_x(R^3\times [0,\infty))$. The definition above is a modification of the definition found in \cite{LemRie} by adding the decay condition (\ref{decay}) to ensure uniqueness. An alternative definition, where (\ref{decay}) is replaced  by a condition on the pressure, can be found in \cite{KiSe}. The existence of Leray solutions for very general inital data is proved in \cite{LemRie}. In our situation with initial data $u_0$ in $L^3$ we can follow \cite{Cald, RuSv} or see section 4 below. We note that condition  (\ref{decay}) allows us to calculate $p$ in the following way: $\forall B_r(x_0)\times (0,t_{\ast})\subseteq R^3\times(0,\infty)$, take a smooth cutoff function $\phi$ with $\phi|_{B_{2r}(x_0)}=1$,  then there exists a function $p(t)$ depending only on $x_0,r,t,\phi$ (we suppress the dependence on $x_0,r,\phi$ in our notation) such that for $(x,t)\in B_r(x_0)\times(0,t_{\ast})$
\begin{equation}\label{eq:no3}
p(x,t)=-\Delta^{-1}\mbox{div~div}\left(u\otimes u\phi\right)-\int_{R^3}\left(k(x-y)-k(x_0-y)\right)u\otimes u(y,t)\left(1-\phi(y)\right)dy+p(t),
\end{equation}
where $k(x)$ is the kernel of $\Delta^{-1}\mbox{div~div}$. \\
The right hand side is well defined since $u$ satisfies the estimates in i) and 
\begin{equation}
|k(x-y)-k(x_0-y)|=O\left(\frac{1}{|x_0-y|^4}\right) {\rm ~~as~~} |y|\to \infty. 
\end{equation}
The situation is similar to extending the domain of singular integrals to bounded functions, see for example \cite{LemRie} and \cite{Stein}.
There are many other possibilities in choosing  a  decay condition (such as by imposing conditions on the pressure used in \cite{KiSe}, already mentioned above). Condition (\ref{decay}) works well for our purposes here.  It should be noted that  some decay of $u$ at spatial infinity is needed if we wish  $p$ to be given as in (\ref{eq:no3}). For instance, as observed by many authors, if we take $u(x,t)=f(t)$, $p(x,t)=-f'(t)\cdot x$, then $(u,p)$ verifies the conditions for the definition of Leray solution except the decay requirement. In this case, $p$ is not given as the above formula though the right hand side is still well defined.\\
We will use the following version of the local energy estimates due to  Lemari\'{e}-Rieusset  \cite{LemRie}:
\begin{lemma}{(A priori estimate for Leray solution)}\label{lem:no1}\\
 Let $\alpha=\sup_{x_0\in R^3}\int_{B_R(x_0)}\frac{|u_0|^2}{2}(x)dx<\infty$ for some $R>0$ and let $u$ be a Leray solution with initial data $u_0$. Then for $\lambda$ satisfying $0<\lambda\leq\epsilon_0\min\{\alpha^{-2}R^2,1\}$ with some small absolute number $\epsilon_0>0$, we have
\begin{equation}\label{eq:no4}
 {\rm ess}\sup_{0\leq t\leq \lambda R^2}\sup_{x_0\in R^3}\int_{B_R(x_0)}\frac{|u|^2}{2}(x,t)dx+\sup_{x_0\in R^3}\int_0^{\lambda R^2}\int_{B_R(x_0)}|\nabla u|^2(x,t)dxdt\leq C\alpha.
\end{equation}

\end{lemma}

\noindent
\textbf{Proof}: Since $u$ is suitable and $u(t,\cdot)$ converges to $u_0$ locally in $L^2$ as $t\to 0+$, we obtain by local energy estimate:\\
\begin{eqnarray*}
 &&\int_{R^3}\frac{|u|^2}{2}(x,t)\phi(x-x_0)dx+\int_0^t\int_{R^3}|\nabla u|^2\phi(x-x_0)dxds\\
&&\leq \int_{R^3}\frac{|u|^2}{2}(x,0)\phi(x-x_0)dx+\int_0^t\int_{R^3}\frac{|u|^2}{2}\Delta \phi(x-x_0)dxds\\
&&~~~~+\int_0^t\int_{R^3}\frac{|u|^2}{2}u\cdot\nabla \phi(x-x_0)+pu\cdot \nabla\phi(x-x_0)dxds,
\end{eqnarray*}
for a.e. $t>0$, where $\phi$ is a nonnegative smooth cutoff function with $\phi=1$ in $B_R(0)$, ${\rm supp}~~\phi\Subset B_{2R}(0)$ and $|\nabla \phi|\leq \frac{C}{R}$. For $\lambda<1$, denote
\begin{equation}
A(\lambda):={\rm ess}\sup_{0\leq t\leq \lambda R^2}\sup_{x_0\in R^3}\int_{R^3}\frac{|u|^2}{2}(x,t)\phi(x-x_0)dx+\sup_{x_0\in R^3}\int_0^{\lambda R^2}\int_{R^3}|\nabla u|^2(x,t)\phi(x-x_0)dxdt\mbox{.}
\end{equation}
Sobolev embedding theorem gives
\begin{equation*}
\sup_{x_0\in R^3}\int_0^{\lambda R^2}\int_{B_{2R}(x_0)}|u|^3(x,t)dxdt\leq C A(\lambda)^{3/2}R^{1/2}\lambda^{1/4}\quad \mbox{if $\lambda\leq 1$},
\end{equation*}
 for some absolute number $C>0$. We apply formula (\ref{eq:no3}) to $p$ in $B_{2R}(x_0)\times (0,\lambda R^2)$:
\begin{equation*}
 p(x,t)=-\Delta^{-1}\mbox{div~~div} \left(u\otimes u\psi\right)-\int_{R^3}\left(k(x-y)-k(x_0-y)\right)\left(u\otimes u(y,t)(1-\psi(y))\right)dy+p(t),
\end{equation*}
where $\psi$ is a smooth cutoff function with $\psi|_{B_{4R}(x_0)}=1$, $0\leq \psi\leq 1$, $\psi$ vanishes outside $B_{8R}(x_0)$ and $|\nabla \psi|\leq \frac{C}{R}$.
Then by elliptic estimates and 
\begin{equation}
|k(x-y)-k(x_0-y)|\leq \frac{CR}{|x_0-y|^4} {\rm ~~for~~} |x_0-y|\ge 4R, |x-x_0|\leq 2R, 
\end{equation}
we easily obtain
\begin{eqnarray*}
&& \|p(x,t)-p(t)\|_{L^{3/2}(B_{2R}(x_0)\times(0,\lambda R^2))}\\
&&\leq C \left(\|u\|_{L^3(B_{8R}(x_0)\times (0,\lambda R^2))}^2+\|R^{-3}A(\lambda)\|_{L^{3/2}(B_{2R}(x_0)\times (0,\lambda R^2))}\right)\\
&&\leq C \lambda^{1/6}A(\lambda)R^{1/3},{\rm ~~for~~} \lambda\leq 1.
\end{eqnarray*}
Thus we obtain from the local energy inequality for $\lambda\leq 1$ and a.\ e.\ $t\leq \lambda R^2$:
\begin{equation*}
\int_{R^3}\frac{|u|^2}{2}(x,t)\phi(x-x_0)dx+\int_0^t\int_{R^3}|\nabla u|^2\phi(x-x_0)dxds\leq \alpha+C\lambda A(\lambda)+CA(\lambda)^{3/2}\lambda^{1/4}R^{-1/2}.
\end{equation*}
Taking sup over $x_0\in R^3$ and $t\leq\lambda R^2$, we get 
\begin{equation}
A(\lambda)\leq \alpha+CA(\lambda)\lambda+CA(\lambda)^{3/2}\lambda^{1/4}R^{-1/2}. 
\end{equation}
Note that $A(\lambda)$ is a priori bounded which is critical in our lemma. Also, we note that $A(\lambda)$ is non-decreasing in $\lambda$ and from (\ref{eq:no2}) it is not hard to see that $A(\lambda)$ is continuous in $\lambda$. (We note that this conclusion does not imply the continuity of the map $t\to \int_{R^3} |u(x,t)|^2\psi(x)dx$, which is unclear, in general.) From the above estimate for $A(\lambda)$, the lemma follows easily by the usual ``continuation in $\lambda$" argument when $\epsilon_0$ is chosen sufficiently small. Note that from the formula (\ref{eq:no3}) and the a priori estimate of $u$, we get the following estimate for $p$ which will be useful: 
\begin{equation}
\sup_{x_0\in R^3}\int_0^{\lambda R^2}\int_{B_R(x_0)}|p-p(t)|^{3/2}dxdt\leq C \alpha^{3/2}R^{1/2}.
\end{equation}

\noindent
\textbf{Remarks:} In the above estimate on $p$, more precisely, $p(t)=p_{x_0,R}(t)$. That is, we need to choose some appropriate constants $p_{x_0,R}(t)$ to satisfy the inequality. The point here is that such constants depending on $x_0,R,t$ exist. This remark is effective throughout the paper. \\

We have the following simple corollary that will be useful below.
\begin{corollary}\label{cor:no1}
 Let $u$ be a Leray solution with initial data $u_0\in L^3(R^3)$. Let $p$ be the associated pressure. Then for $\forall r>0$,
\begin{eqnarray}
&&\int_0^{r^2}\int_{B_r(x_0)}|\nabla u|^2dxds+{\rm ess}\sup_{0\leq t\leq r^2}\int_{B_r(x_0)}\frac{|u|^2(x,t)}{2}dx\nonumber\\
&&\leq \frac{C\|u_0\|_{L^3(R^3)}^2r}{\sqrt{\epsilon_0\min\{\|u_0\|_{L^3(R^3)}^{-4},1\}}}{\rm ,~~~and}\\
&&\int_0^{r^2}\int_{B_r(x_0)}|p-p(t)|^{3/2}dxds\leq \frac{C\|u_0\|_{L^3(R^3)}^3r^2}{\epsilon_0\min\{\|u_0\|_{L^3(R^3)}^{-4},1\}},
\end{eqnarray}
for any $x_0\in R^3$.
\end{corollary}

\noindent
\textbf{Proof:} For each $r>0$, let $R=\frac{r}{\sqrt{\epsilon_0\min\{\|u_0\|_{L^3(R^3)}^{-4},1\}}}>r$. We shall apply Lemma \ref{lem:no1} with this $R$. We have:
\begin{eqnarray*}
\alpha&=&\sup_{x_0\in R^3}\int_{B_R(x_0)}|u_0|^2dx\leq \sup_{x_0\in R^3}\left(\int_{B_R(x_0)}|u_0|^3dx\right)^{2/3}R\\
&\leq&\|u_0\|_{L^3(R^3)}^2R.
\end{eqnarray*}
Thus we can choose $\lambda=\epsilon_0\min\{\|u_0\|_{L^3(R^3)}^{-4},1\}\leq \epsilon_0\min\{\alpha^{-2}R^2,1\}$. Note that by our choice of $R$ and $\lambda$, $\lambda R^2=r^2$, therefore from Lemma \ref{lem:no1}, the lemma follows.\\

We will prove the following uniqueness result, which is a variant of classical uniqueness results (such as \cite{Serrin2, LemRie}).
The situation we consider here is on one hand slightly more general than the one considered in \cite{Serrin2} and, on the other hand, simpler than the one considered in \cite{LemRie}, allowing for a simple proof (still based on similar ideas).
 
\begin{lemma}\label{lem:no2}
Suppose $(u,p)$, $(v,q)$ are two Leray solutions with the same initial data $u_0$, and $v\in L^5(R^3\times [0,T))$ for any $T<\infty$. Then $u=v$ in $R^3\times (0,\infty)$ almost everywhere.
 \end{lemma}

\noindent
\textbf{Proof:} It suffices to prove $u=v$ almost everywhere in $R^3\times [0,T)$ for any $T<\infty$. Since $u,v$ satisfy NSE, we can change $u,v$ in a set of measure zero, such that for any $B_R(0)$, $t\rightarrow u(t)$, $t\rightarrow v(t)$ are weakly continuous in $L^2(B_R(0))$ for $t\in [0,T)$ (in fact, $t\to v(\cdot,t)$ is strongly continuous in $L^2(B_R(0))$ by the Serrin-Ladyzhenskaya-Prodi regularity criteria \cite{Ladyzhenskaya}\cite{Serrin}\cite{Prodi}).\\
We first observe that for any smooth compactly supported $\psi\ge 0$, $\int_{R^3}u(x,t)\cdot v(x,t)\psi(x)dx$ as a function in $t$ is in $W^{1,1}([0,T))$, and
\begin{eqnarray}\label{eq:no5}
 &&\frac{d}{dt}\int_{R^3}u(x,t)\cdot v(x,t)\psi(x)dx \nonumber\\
&&=\int_{R^3}-2(\nabla u:\nabla v)(x,t)\psi(x)-(\nabla u:v\otimes  \nabla \phi)(x,t)-(u\cdot\nabla u) v\psi dx\nonumber\\
&&~~~~+\int_{R^3}pv\cdot \nabla \psi -(\nabla v:u\otimes\nabla\psi)(x,t)-(v\cdot\nabla v) u \psi+qu\cdot \nabla\psi dx,
\end{eqnarray}
in the sense of distributions on $(0,T)$.\\
The right hand side is easily checked to be integrable in $t$ using $p,q\in L^{3/2}_{loc}$, $\nabla u, ~~\nabla v\in L^2_{loc}$,
$u\in L^{10/3}_{loc}$, $v\in L^5_{loc}$. Then the proof follows by a usual mollification procedure.\\
Next, from local energy inequality and the above identity integrated in time, we have for any smooth compactly supported $\psi\ge 0$ and a.\ e.\ $t$
\begin{eqnarray*}
 &&\int_{R^3}\frac{|u-v|^2}{2}(x,t)\psi(x)dx\\
&&\leq \int_0^t\int_{R^3}\frac{|u|^2}{2}\Delta \psi+\frac{|u|^2}{2}u\cdot \nabla \psi+pu\cdot \nabla\psi dxds\\
&&~~~~+\int_0^t\int_{R^3}\frac{|v|^2}{2}\Delta \psi+\frac{|v|^2}{2}v\cdot \nabla \psi+qv\cdot \nabla\psi dxds\\
&&~~~~+\int_0^t\int_{R^3}\left(-|\nabla u|^2\psi -|\nabla v|^2\psi +2\nabla u:\nabla v \psi\right) dxds\\
&&~~~~+\int_0^t\int_{R^3}(\nabla u:v\otimes\nabla \psi)+(u\cdot\nabla u)\cdot v\psi-pv\cdot\nabla\psi dxds\\
&&~~~~+\int_0^t\int_{R^3}(\nabla v:u\otimes\nabla \psi)+(v\cdot\nabla v) u\psi-qu\cdot\nabla\psi dxds\mbox{.}
\end{eqnarray*}
Note that we have the following relatitions:
\begin{eqnarray*}
 &&\int_0^t\int_{R^3}(\nabla u:v\otimes\nabla \psi)+(\nabla v:u\otimes\nabla \psi) dxds\\
&&=\int_0^t\int_{R^3}\nabla(u\cdot v)\nabla\psi dxds=-\int_0^t\int_{R^3}u\cdot v\Delta \psi dxds;
\end{eqnarray*}
and,
\begin{eqnarray*}
&&\int_0^t\int_{R^3}(u\cdot\nabla u)\cdot v\psi+(v\cdot\nabla v)\cdot u\psi dxds\\
&&=\int_0^t\int_{R^3}\left(u\cdot\nabla(u-v)\right)v\psi+\left(v\cdot\nabla(v-u)\right)u\psi dxds\\
&&~~~~-\int_0^t\int_{R^3}\left(\frac{|v|^2}{2}u\cdot\nabla\psi+\frac{|u|^2}{2}v\cdot\nabla\psi\right) dxds\\
&&=\int_0^t\int_{R^3}\left((u-v)\cdot\nabla(u-v)\right)v\psi +\left(v\cdot\nabla(v-u)\right)(v-u)\psi dxds\\
&&~~~~-\int_0^t\int_{R^3}\left(\frac{|v|^2}{2}u\cdot\nabla\psi+\frac{|u|^2}{2}v\cdot\nabla\psi\right) dxds.
\end{eqnarray*}
Various integration by parts in the above hold since $|v||\nabla u||u|\in L^1_{loc}$. Here $v\in L^5_{loc}$ plays a critical role.
Thus, summarizing the above, we get:
\begin{eqnarray*}
 &&\int_{R^3}\frac{|u-v|^2}{2}(x,t)\psi(x)dx\\
&&\leq -\int_0^t\int_{R^3}|\nabla (u-v)|^2\psi(x)+\frac{|u-v|^2}{2}\Delta \psi dxds+\int_0^t\int_{R^3}(p-q)(u-v)\cdot\nabla\psi dxds\\
&&~~~~+\int_0^t\int_{R^3}\left((u-v)\cdot\nabla(u-v)\right) v\psi dxds+\int_0^t\int_{R^3}\left(v\cdot\nabla(u-v)\right)(v-u)\psi dxds\\
&&~~~~+\int_0^t\int_{R^3}(\frac{|u|^2}{2}-\frac{|v|^2}{2})(u-v)\cdot\nabla \psi dxds\\
&&\leq -\int_0^t\int_{R^3}|\nabla (u-v)|^2\psi(x)+\frac{|u-v|^2}{2}\Delta \psi dxds+\int_0^t\int_{R^3}(p-q)(u-v)\cdot\nabla\psi dxds\\
&&~~~~+\int_0^t\int_{R^3}\left((u-v)\cdot\nabla(u-v)\right)\psi dxds+\int_0^t\int_{R^3}\left(v\cdot\nabla(u-v)\right)(v-u)\psi dxds\\
&&~~~~+\int_0^t\int_{R^3}\frac{|u-v|^3}{2}|\nabla \psi|dxds+\int_0^t\int_{R^3}|v||u-v|^2|\nabla\psi|dxds\mbox{.}
\end{eqnarray*}
Denote 
\begin{equation}
e(t)={\rm ess}\sup_{0\leq s\leq t,~x_0\in R^3}\int_{B_1(x_0)}\frac{|u-v|^2}{2}(x,s)dx+\sup_{x_0\in R^3}\int_0^t\int_{B_1(x_0)}|\nabla (u-v)|^2dxds.
\end{equation}
By multiplicative inequalities, we have
\begin{equation}
 \sup_{x_0\in R^3}\|u-v\|_{L^{10/3}(B_1(x_0)\times(0,t))}\leq C e(t)^{1/2},{\rm ~~~~for~~}t\leq 1.
\end{equation}
For any $x_0$, let $\psi$ be a standard smooth cutoff function with $\psi|_{B_1(x_0)}=1$, $\psi$ vanishes outside $B_2(x_0)$. From the formula for $p,q$, we see that for a.e $x\in B_2(x_0)$, $s\in(0,t)$,
\begin{eqnarray*}
(p-q)(x,s)&=&p(s)-q(s)-\Delta^{-1}\mbox{div~~div}\left(u\otimes (u-v)\phi\right)-\Delta^{-1}\mbox{div~~div}\left((u-v)\otimes v\phi\right)\\
&&~~~~-\int_{R^3}\left(k(x-y)-k(x_0-y)\right)\left(u\otimes(u-v)+(u-v)\otimes v\right)(1-\phi)dy,
\end{eqnarray*}
for some standard cutoff function $\phi$ with $\phi|_{B_4(x_0)}=1$ and $\phi$ vanishing outside $B_8(x_0)$.\\
Thus we obtain
\begin{eqnarray*}
 &&\|p-q-\left(p(s)-q(s)\right)\|_{L^{3/2}(B_2(x_0)\times(0,t))}\\
&&\leq C \left(\|u\|_{L^3(B_8(x_0)\times(0,t))}+\|v\|_{L^3(B_8(x_0)\times(0,t))}\right)\|u-v\|_{L^3(B_8(x_0)\times(0,t))}\\
&&~~~~+C\left(\int_0^t\sup_{y\in R^3}\|u-v\|_{L^2(B_1(y))}^{3/2}(s)ds\right)^{2/3}\sup_{y\in R^3,0\leq s\leq t}\left(\|u\|_{L^2(B_1(y))}(s)+\|v\|_{L^2(B_1(y))}(s)\right)\\
&&\leq C(u,v,T)e(t)^{1/2}t^{1/15}\quad \mbox{for $t\leq 1$.}
\end{eqnarray*}
Since $e(t)\leq C(u,v)$ we have
\begin{eqnarray*}
 &&\int_{R^3}\frac{|u-v|^2}{2}(x,t)\psi(x)dx+\int_0^t\int_{R^3}|\nabla(u-v)|^2\psi(x)dx\\
&&\leq C(u,v,T)\left(e(t)t^{1/10}+\|v\|_{L^5(R^3\times(0,t))}e(t)\right)\mbox{.}
\end{eqnarray*}
Since we have the freedom to choose $\psi$, by varying the support, and taking supremum, we obtain
\begin{equation*}
 e(t)\leq C(u,v,T)e(t)\left(t^{1/10}+\|v\|_{L^5(R^3\times (0,t))}\right)\quad\mbox{for $t\leq 1$.}
\end{equation*}
This forces $e(t)=0$ for $t<T_{\ast}=T_{\ast}(u,v,T)$, for some sufficiently small $T_{\ast}>0$. After applying this result several times, we see the proof of the lemma is complete.

\medskip
\noindent
\textbf{Remarks:} It is clear from the proof that we only need
\begin{equation} 
\sup_{x_0\in R^3}\int_0^T\int_{B_1(x_0)}|v|^5(x,t)dxdt<\infty \quad {\rm and~}\lim_{x_0\to\infty}\int_0^T\int_{B_1(x_0)}|v|^5(x,t)dxdt=0
\end{equation}
to guarantee uniqueness on $R^3\times [0,T)$.\\

The following version of $\epsilon$-regularity criteria of Caffarelli-Kohn-Nirenberg will be important for us in the sequel:
\begin{lemma}\label{lm:regularity}
Let $(u,p)$ be a suitable weak solution to NSE in $Q_1:=B_1(0)\times(-1,0)$ with $u\in L^{\infty}_tL^2_x(Q_1)\cap L^2_t\dot{H}^1(Q_1)$ and $p\in L^{3/2}(Q_1)$, in the sense that $(u,p)$ verifies NSE as distributions and they satisfy local energy inequality. Then there exists an absolute constant $\epsilon_0>0$, with the following property:\\
if $(\int_{Q_1}|u|^3dxdt)^{1/3}+(\int_{Q_1}|p|^{3/2}dxdt)^{2/3}\leq \epsilon_0$, then $\|\nabla^k u\|_{L^{\infty}(Q_{1/2})}\leq C_k$ for some constants $C_k$, $k=0,1,\dots$
\end{lemma}
A sketch of a short proof can be found for example in \cite{Lin}, a detailed one in \cite{LaSe}.\\

We recall the following lemmas proved in \cite{RuSv}:
\begin{lemma}{(compactness)}\label{lm:compactness}
Let $(u^k,p^k)$, $k=1,2,\dots$ be a sequence of suitable weak solutions such that $u^k$ are uniformly bounded in the energy space $L^{\infty}_tL^2_x\cap L^2_t\dot{H}^{1}_x$ on compact subsets of open set $\mathcal{O}\subset R^3\times R$ and $p^k$ are uniformly bounded in $L^{3/2}_tL^{3/2}_x$ on compact subsets of $\mathcal{O}$. Then the sequence $u^k$ is compact in $L^3_tL^3_x$ on compact subsets of $\mathcal{O}$. Moreover, if $u^k\to u$ in $L^3_tL^3_x$ on compact subsets of $\mathcal{O}$ and $p^k\rightharpoonup p$ in $L^{3/2}_tL^{3/2}_x$ on compact subsets of $\mathcal{O}$, then $(u,p)$ is again a suitable weak solution.
\end{lemma}

\begin{lemma}{(Stability of Singularities)}\label{lm:stability}
 In the situation of lemma \ref{lm:compactness}, assume that $z^k\in \mathcal{O}$ are singular points of $(u^k,p^k)$, $k=1,2,\dots$, and that $z^k\to z_0\in \mathcal{O}$. Then $z_0$ is a singular point of $(u,p)$.
\end{lemma}
We refer readers to \cite{RuSv} for the proofs, and here we only recall that a point $z_0$ is called a singular point of a suitable weak solution $u$ to NSE if $u$ is not bounded in any neighborhood of $z_0$.
\end{section}

\begin{section}{Mild solutions with initial data in $L^3(R^3)$}
In this section we collect some well known results about mild solutions, and introduce some splitting arguments which are useful in the proof of our main result.\\
One can rewrite NSE as an integral equation:
\begin{equation}\label{eq:integralequation}
 u(\cdot,t)=e^{\Delta t}u_0-\int_0^te^{\Delta (t-s)}P\mbox{div}~u\otimes u(\cdot,s)ds,
\end{equation}
where $P$ is the Helmholtz projection operator. It is well known that (\ref{eq:integralequation}) has a global solution if initial data is small in $L^3$, and for arbitrary initial data in $L^3$ a unique local in time solution $u\in C([0,T_{\ast}),L^3(R^3))$ with a number of additional properties such as $u\in L^5(R^3\times (0,T))$ for $T<T_{\ast}$, here $T_{\ast}$ denotes the maximal existence time. See \cite{Kato, KoTa} and references therein for more details. Take any $v\in \mathcal{N}(u_0)$, by the uniqueness Lemma \ref{lem:no2} and the remark below it, we know $v=u$ on $R^3\times [0,T_{\ast})$. Thus when we consider properties of solutions only on $R^3\times [0,T_{\ast})$, there is no confusion to assume $u$ has been properly extended to $R^3\times[0,\infty)$ as a Leray solution. We will make use of this observation below.\\

A priori there could be a number of reasons why $T_{\ast}$ could be finite, we first show that it can only be due to the formation of a `singular point':
\begin{lemma}\label{lem:sing}
Let $u\in C([0,T_{\ast}),L^3(R^3))$ be the mild solution to NSE with initial data $u_0$ and $T_{\ast}$ is the maximal existence time. Suppose $T_{\ast}<\infty$, then there exists $z_0=(x_0,T_{\ast})$ such that $\forall r>0$, ${\rm ess}\sup_{Q_r(z_0)}~|u|=+\infty$.
\end{lemma}

\noindent
\textbf{Proof:} We write $u_0=a+b$, with $\|a\|_{L^3(R^3)}\leq \epsilon$, $\|b\|_{L^2(R^3)}<\infty$, where $\epsilon$ is a sufficiently small number to be chosen later. There are a number of ways in which one can perform such a decomposition. One can for example take 
\begin{equation}
a=P(u_0I_{|u_0|<\lambda}), ~~~b=P(u_0I_{|u_0|\ge \lambda}), 
\end{equation}
and take $\lambda$ sufficiently small. If we choose $\epsilon$ small enough, we can apply global existence result for small data for NSE and get a global mild solution $v$ with initial data $a$. Then $w=u-v$ satisfies
\begin{eqnarray}
 \left.\begin{array}{rl}
        \partial_tw-\Delta w+v\cdot\nabla w+w\cdot\nabla v+w\cdot\nabla w+\nabla q&=0\\
        \mbox{div}~~w&=0
       \end{array}\right\} &&\mbox{for $(x,t)\in R^3\times (0,T_{\ast})$,}\label{eq:no7}\\
   w(\cdot,0)=b &&\mbox{in $R^3$.}\label{eq:no11}
\end{eqnarray}
We claim the following estimate:
\begin{equation}\label{estimate}
{\rm ess}\sup_{0\leq t\leq T_{\ast}}\|w(\cdot,t)\|_{L^2(R^3)}^2+\int_0^{T_{\ast}}\int_{R^3}|\nabla w|^2(x,t)dxdt\leq C(\|b\|_{L^2(R^3)},\|v\|_{L^5(R^3\times (0,T_{\ast}))},T_{\ast})\mbox{.}
\end{equation}
There are two ways in which we can prove this claim. Since this type of argument will be used several times we provide them both here.\\
In the first approach, note that regular solution $\tilde{w}$ (by `regular' we mean $\tilde{w}$ is smooth with sufficient decay) to equations (\ref{eq:no7}) (\ref{eq:no11}) satisfies the following, by a simple integration by parts:
\begin{eqnarray*}
 &&\frac{d}{dt}\int_{R^3}\frac{|\tilde{w}|^2}{2}(x,t)dx+\int_{R^3}|\nabla \tilde{w}|^2(x,t)dx\\
&&\leq \int_{R^3}(\tilde{w}\cdot\nabla \tilde{w})v(x,t)dx\\
&&\leq\|v(\cdot,t)\|_{L^5_x}\|\tilde{w}(\cdot,t)\|_{L^{10/3}_x}\|\nabla\tilde{w}(\cdot,t)\|_{L^2_x}\\
&&\leq C\|v(\cdot,t)\|_{L^5_x}\|\tilde{w}(\cdot,t)\|_{L^2_x}^{2/5}\|\nabla\tilde{w}(\cdot,t)\|_{L^2_x}^{8/5}\\
&&\leq \frac{1}{2}\int_{R^3}|\nabla \tilde{w}|^2(x,t)dx+C\|v(\cdot,t)\|_{L^5}^5\int_{R^3}|\tilde{w}(\cdot,t)|^2dx\mbox{.}
\end{eqnarray*}
In the above, we have used H\"{o}lder inequality, interpolation inequality 
\begin{equation*}
\|\tilde{w}\|_{L^{10/3}_x}\leq \|\tilde{w}\|_{L^2_x}^{2/5}\|\tilde{w}\|_{L^6_x}^{3/5},
\end{equation*}
and Sobolev embedding $\dot{H}^1\hookrightarrow L^6$ in $R^3$.\\
Thus we obtain:
\begin{equation*}
 \frac{d}{dt}\int_{R^3}\frac{|\tilde{w}|^2}{2}(x,t)dx\leq C\|v(\cdot,t)\|_{L^5_x}^5\int_{R^3}|\tilde{w}|^2dx\mbox{.}
\end{equation*}
Since $\int_0^{T_{\ast}}\|v(\cdot,t)\|_{L^5_x}^5dt$ is bounded, we get from Gronwall's inequality:
\begin{equation}\label{eq:no12}
 \sup_{0\leq t\leq T_{\ast}}\|\tilde{w}(\cdot,t)\|_{L^2(R^3)}^2+\int_0^{T_{\ast}}\int_{R^3}|\nabla\tilde{w}|^2(x,t)dxdt\leq C(\|b\|_{L^2(R^3)},\|v\|_{L^5(R^3\times(0,T_{\ast}))},T_{\ast})\mbox{.}
\end{equation}
With this a priori estimate at hand, we can then follow Leray's arguments in constructing global weak solutions and obtain a weak solution $\tilde{w}$ to equations (\ref{eq:no7}) (\ref{eq:no11}) and $\tilde{w}$ satisfies the energy inequality (\ref{eq:no12}). We can also require $\tilde{w}$ to satisfy the appropriate local energy inequality (see (\ref{eq:no14}) below). Since $v+\tilde{w}$ is also a Leray solution with initial data $u_0$, by uniqueness result of lemma (\ref{lem:no2}) and the remark below it  we must have $w=\tilde{w}$ on $R^3\times[0,T_{\ast})$ and thus $w$ satisfy (\ref{eq:no12}).\\
Alternatively one can also derive the inequality directly from the equation as follows.\\
$w$ clearly satisfies the following integral equation:
\begin{equation}\label{eq:no13}
w(\cdot,t)=e^{\Delta t}w_0-\int_0^te^{\Delta (t-s)}P\mbox{div~}\left(v\otimes w(\cdot,s)+w\otimes v(\cdot,s)+w\otimes w(\cdot,s)\right)ds,
\end{equation}
where $w_0=b\in L^2(R^3)$, $t<T_{\ast}$.\\
Since both $u,v\in C([0,T_{\ast}),L^3(R^3))$ we conclude $w\in C([0,T_{\ast}),L^3(R^3))$. From this fact and the integral equation (\ref{eq:no13}) together with $b\in L^2(R^3)$ we obtain $w\in C([0,T_{\ast}),L^2(R^3))$ from known estimates of the integral equation. Clearly, $w$ satisfies the following local energy inequality (where $q$ is the associated pressure for $w$):
\begin{eqnarray}
 &&\int_{R^3}\frac{|w|^2}{2}(x,t)\phi(x)dx+\int_0^t\int_{R^3}|\nabla w|^2\phi(x)dx\nonumber\\
&&\leq \int_{R^3}\frac{|w|^2}{2}(x,0)\phi(x)dx+\int_0^t\int_{R^3}\frac{|w|^2}{2}\Delta \phi+qw\cdot\nabla\phi dxds\nonumber\\
&&~~~~+\int_0^t\int_{R^3}\frac{|w|^2}{2}(w+v)\cdot\nabla \phi-(w\cdot\nabla v)w\phi dxds, \label{eq:no14}
\end{eqnarray}
where $\phi$ is a nonnegative smooth cutoff function.
By local theory of mild solutions (see \cite{Kato}), we know 
\begin{equation}
\sup_{0\leq t\leq T_{\ast}}\sqrt{t}\|\nabla v(\cdot,t)\|_{L^3_x}<\infty. 
\end{equation}
Thus 
\begin{equation}
\int_0^T\int_{R^3}|w||\nabla v||w|dxds\leq \int_0^T\|\nabla v(\cdot,t)\|_{L^3_x}\|w(\cdot,t)\|_{L^3}^2dx<\infty,
\end{equation}
for any $T<T_{\ast}$. This, together with $q\in L^{\infty}((0,T),L^{3/2}(R^3)),~w\in C([0,T),L^3\cap L^2(R^3))$ for any $T<T_{\ast}$, implies that we can take the
cutoff function in the local energy inequality (\ref{eq:no14}) to be $\phi(\frac{x}{R})$ with $\phi|_{B_1}\equiv 1$, and send $R\to \infty$. We obtain
\begin{eqnarray*}
 &&\int_{R^3}\frac{|w|^2}{2}(x,t)dx+\int_0^t\int_{R^3}|\nabla w|^2dxds\\
&&\leq \int_{R^3}\frac{|w|^2}{2}(x,0)dx-\int_0^t\int_{R^3}(w\cdot\nabla v)wdxds\\
&&=\int_{R^3}\frac{|w|^2}{2}(x,0)dx+\int_0^t\int_{R^3}(w\cdot\nabla w)vdxds.
\end{eqnarray*}
The last identity holds since the first inequality implies 
\begin{equation*}
\int_0^t\int_{R^3}|\nabla w|^2dxds<\infty, 
\end{equation*}
which is the key. Then we can proceed as in the first approach to finish the proof of the claim.\\
From Corollary \ref{cor:no1}, we know
\begin{equation}
\sup_{x_0\in R^3,0\leq t\leq T_{\ast}}\int_{B_1(x_0)}|v|^2(x,t)dx+\sup_{x_0\in R^3}\int_0^{T_{\ast}}\int_{B_1(x_0)}|\nabla v|^2+|\tilde{p}-\tilde{p}(t)|^{3/2}dxdt\leq C(\epsilon)C(T_{\ast}),
\end{equation}
where $\tilde{p}$ is the associated pressure for $v$. Observe from equation (\ref{eq:no7}) and estimate (\ref{estimate}) that $q\in L^{3/2}(R^3\times(0,T_{\ast}))$ and $w\in L^3(R^3\times(0,T_{\ast}))$. Since $C(\epsilon)\to 0$ as $\epsilon\to 0$, we can choose $\epsilon$ sufficiently small and find $R>0$ sufficiently large, such that
for $|x_0|>R$, we can apply the $\epsilon$-regularity criteria to $u=v+w$ in $Q_{\min(\sqrt{T_{\ast}/2},1)}(x_0,T_{\ast})$. Thus we are able to conclude the following:

\medskip
\noindent
{\sl \,\, there exists a compact set $K\subseteq R^3$ such that $u$ is bounded in $(R^3\backslash K)\times [T_{\ast}/2,T_{\ast}]$.}

\medskip
\noindent
Now suppose the lemma is not true. Then $u$ is bounded in a neighborhood of any point in $R^3\times[T_{\ast}/2,T_{\ast}]$. Thus by the compactness of $K$, we see $u$ is bounded in $R^3\times[T_{\ast}/2,T_{\ast}]$. Since $v$ is bounded in $R^3\times[T_{\ast}/2,T_{\ast}]$, we conclude $w$ is also bounded in $R^3\times[T_{\ast}/2,T_{\ast}]$. This means $u(\cdot,t)$ is bounded in $L^3$ as well as $L^{\infty}$, as $t$ approaches $T_{\ast}$. Then local existence theory for NSE tells us we can continue $u$ beyond $T_{\ast}$, a contradiction.

\begin{remark}\label{remarks}
The splitting argument in the above proof is very useful in obtaining estimates for $u$ even when we approach the blow up time. Let $u\in C([0,1),L^3(R^3))$ be a mild solution with initial data $u_0$. Still take the above decomposition $u_0=a+b$ with $a,b$ defined as above. Denote $\alpha=\|u_0\|_{L^3}$. We immediately see 
\begin{equation}
\|a\|_{L^6(R^3)}\leq C\lambda^{1/2}\alpha^{1/2}, ~~\|b\|_{L^2(R^3)}\leq C\lambda^{-1/2}\alpha^{3/2}. 
\end{equation}
If we choose $\lambda=\lambda(\alpha)$ so small such that $\|a\|_{L^6(R^3)}$ is smaller than some absolute number, fix such $\lambda=\lambda(\alpha)$, then we can conclude from local existence theory of NSE the existence of a mild solution $v$ to NSE with initial data $a$ in $R^3\times[0,2)$ with 
\begin{equation}
\,\,{\rm ess~sup~~}_{0\leq t\leq 3/2}\|v(\cdot,t)\|_{L^6}\leq C. 
\end{equation}
Then $w$ as defined above satisfy (\ref{eq:no7}) with initial data $b$. Now we can use the usual energy estimate as in the above lemma to bound 
\begin{equation}
\sup_{0\leq t\leq 1}\|w(\cdot,t)\|_{L^2(R^3)}^2+\int_0^1\int_{R^3}|\nabla w|^2dxdt\leq C(\alpha). 
\end{equation}
Thus, $u=v+w$ satisfies 
\begin{equation}
\sup_{0\leq t\leq 1}\|u(\cdot,t)\|_{L^2+L^6}\leq C(\alpha). 
\end{equation}
Of course one can do the same thing if the blowup time is $T$, though the estimates will also depend on $T$ then. The point here is this estimate is uniform as long as $L^3$ norm of $u_0$ stays bounded. This observation will be useful later. And in fact, this estimate is stronger than the one in Lemma \ref{lem:no1} for general Leray solutions, since it implies decay of $u$ at spatial infinity while Lemma \ref{lem:no1} does not imply any decay of solutions. The difference is that here we are dealing with an a priori regular solution.
\end{remark}

The following estimate is the main new observation  that enables us to work in $L^3$. 
\begin{lemma}\label{lem:key}
Let $u$ be a Leray solution with divergence free initial data $u_0\in L^3(R^3)$. Then there exists a nonnegative function $h(t)$ depending only on $\|u_0\|_{L^3(R^3)}$, such that $\lim_{t\to 0+}h(t)=0$ and 
\begin{equation}
\|u(\cdot,t)-e^{\Delta t}u_0\|_{L^2(B_1(x_0))}\leq h(t),
\end{equation}
for any $x_0\in R^3$, a.e. $0\leq t<1$.
\end{lemma}

\noindent
\textbf{Proof:} We use a splitting argument which is slightly different from the usual ones. We refer the reader to the paper \cite{GS1} for another example of a splitting argument based on a comparison with the linear equation. Denote $\alpha:=\|u_0\|_{L^3(R^3)}$. For any $\epsilon>0$, we split $u_0=a+b$, with 
\begin{equation}
a=P(u_0I_{|u_0|<M}),~~b=P(u_0I_{|u_0|\ge M}),
\end{equation}
where $M$ is some large number to be chosen later. Clearly, 
\begin{eqnarray*}
&&\|a\|_{L^6(R^3)}\leq C\sqrt{M}\alpha^{1/2},\\
&&\|b\|_{L^2(R^3)}\leq \frac{C}{\sqrt{M}}\|u_0\|_{L^3(R^3)}^{3/2}=\frac{C}{\sqrt{M}}\alpha^{3/2}.
\end{eqnarray*}
Choose $M$ such that 
\begin{equation}
\frac{C}{\sqrt{M}}\alpha^{3/2}<\frac{\epsilon}{100}, 
\end{equation}
and we fix this $M$ from now on, thus $M=M(\epsilon,\alpha)$. From the local existence theory of NSE, we can find $T=T(\epsilon,\alpha)>0$, and a mild solution $v\in C([0,T),L^6(R^3))$ for NSE with initial data $a$, enjoying a number of other properties. Among them in particular, we have
\begin{equation}
\sup_{0\leq t<T}\|v(\cdot,t)\|_{L^6(R^3)}\leq C(\epsilon,\alpha).
\end{equation}
$w=u-v$ satisfies:
\begin{eqnarray}\label{eq:no8}
 \left.\begin{array}{rl}
        \partial_tw-\Delta w+v\cdot\nabla w+w\cdot\nabla v+w\cdot\nabla w+\nabla q&=0\\
        \mbox{div}~~w&=0
       \end{array}\right\} &&\mbox{$(x,t)\in R^3\times (0,T)$,}\\
   w(\cdot,0)=b &&\mbox{in $R^3$.}
\end{eqnarray}
By estimates of $u$ from Corollary \ref{cor:no1} and estimates on $v$, local energy estimates for $w$ (\ref{eq:no14}) and parabolic regularity, one can conclude (we omit the routine calculations):
there exists $T_1(\alpha,\epsilon)>0$ such that
$$\sup_{x_0\in R^3, t\leq T_1(\alpha,\epsilon)}\|v(\cdot,t)-e^{\Delta t}a\|_{L^2(B_1(x_0))}\leq \frac{\epsilon}{10}\mbox{,}$$
$$\sup_{x_0\in R^3, t\leq T_1(\alpha,\epsilon)}\|w(\cdot,t)\|_{L^2(B_1(x_0))}\leq \frac{\epsilon}{2}\mbox{.}$$
Thus, from $u=v+w$, we obtain 
\begin{equation}
\|u(\cdot,t)-e^{\Delta t}u_0\|_{L^2(B_1(x_0))}<\epsilon, 
\end{equation}
for any $x_0\in R^3$, a.e. $t\leq T_1(\alpha,\epsilon)$. From this, the lemma follows easily.

\end{section}

\begin{section}{The main theorem}
For any divergence free $u_0\in L^3(R^3)$, denote $T_{\max}(u_0)$ as the maximal time of existence for the mild solution for NSE starting from $u_0$. Define $\rho_{\max}=\sup\{\rho: T_{\max}(u_0)=\infty ~~\mbox{for every divergence free}~~ u_0\in L^3(R^3) ~~\mbox{with}~~\|u_0\|_{L^3(R^3)}<\rho\}$. Also define $\mathcal{M} :=\{u_0\in L^3(R^3): T_{\max}(u_0)<\infty, \|u_0\|_{L^3(R^3)}=\rho_{\max}\}$.
\begin{theorem}
 Suppose $\rho_{\max}<\infty$. Then $\mathcal{M}$ is nonempty, and moreover, $\mathcal{M}$ is  compact with respect to $L^3$-norm modulo translations and scalings. That is, for any sequence $u_0^k$ in $\mathcal{M}$, there exist $x_k$, $\lambda^k$ such that
$\lambda^ku_0^k(\lambda^k(x-x_k))$ has a convergent subsequence in $L^3(R^3)$.
\end{theorem}

\noindent
\textbf{Proof:} By the definition of $\rho_{\max}$ and the assumption that $\rho_{\max}<\infty$, there exists a sequence of divergence free initial data $u_0^k$ such that $T_{\max}(u_0^k)<\infty$ (thus $\|u_0^k\|_{L^3(R^3)}\ge \rho_{\max}$) and $\|u_0\|_{L^3(R^3)}\to \rho_{\max}$.
By Lemma (\ref{lem:sing}) we know there are singular points for mild solutions $u^k$ corresponding to $u_0^k$. By translations and scalings 
\begin{equation*}
u_0^k\rightarrow \lambda^ku_0^k(\lambda^k(x-x_k)),
\end{equation*}
for some $\lambda^k,x_k$, we can assume the first singularity is at time $1$ and is $(x,t)=(0,1)$. We still denote the sequence as $u_k~(u_0^k~~\mbox{correspondingly})$ after translations and scalings.
By Lemma (\ref{lem:key}), we have
\begin{equation}\label{eq:estimate}
\sup_{x_0\in R^3}\|u^k(\cdot,t)-e^{\Delta t}u_0^k\|_{L^2{(B_1(x_0))}}\leq h(t)\quad\mbox{for $t<1$,}
\end{equation}
for some nonnegative function $h(t)$ with $\lim_{t\to 0+}h(t)=0$.
Note that Corollary \ref{cor:no1} and Remark \ref{remarks} imply uniform boundedness of $u^k$ in local energy norm and 
\begin{equation}\label{eq:decay}
\sup_{0\leq t\leq 1}\|u^k(\cdot,t)\|_{L^2+L^6}\leq C(\rho_{\max}). 
\end{equation}
By compactness as in Lemma \ref{lm:compactness}, and weak continuity in $t$ we can find a subsequence of $u^k$ (which we still denote as $u^k$) and a suitable weak solution $u$ to NSE, such that:\\
$u^k\to u$ in $L^3(B_1(x_0)\times (0,1))$, for all $x_0\in R^3$,\\
$u^k(\cdot,t)\rightharpoonup u(\cdot,t)$ in $L^2(B_1(x_0))$ for every $t\in [0,1]$ and $x_0\in R^3$, and\\ $u_0^k\rightharpoonup u_0$ in $L^3(R^3)$. \\
Moreover, by stability of singularity in Lemma \ref{lm:stability} $(x,t)=(0,1)$ is a singular point of $u$. From estimate (\ref{eq:estimate}) and the weak convergence of $u^k(\cdot,t),~u_0^k$, we get 
\begin{equation}
\sup_{x_0\in R^3}\|u(\cdot,t)-e^{\Delta t}u_0\|_{L^2(B_1(x_0))}\leq h(t). 
\end{equation}
Since $h(t)\to 0$ as $t\to 0$, we see 
\begin{equation*}
u(\cdot,t)\to u_0 {\rm ~~in~~} L^2(B_1(x_0)) {\rm ~~for~~ any~~} x_0.
 \end{equation*}
Furthermore, 
\begin{equation}
\sup_{0\leq t\leq 1}\|u(\cdot,t)\|_{L^2+L^6}\leq C(\rho_{\max}), 
\end{equation}
by the weak convergence of $u^k$ and the estimate (\ref{eq:decay}). Thus $u$ satisfies the decay condition at spatial infinity required in the definition of Leray solutions. Summarzing above, we see $u$ is a Leray solution with initial data $u_0$. By uniqueness result of Lemma \ref{lem:no2} and the remark below it, we see that the mild solution starting with $u_0$ must have a singularity in $R^3\times [0,1]$. Since 
\begin{equation}
\|u_0\|_{L^3(R^3)}\leq\liminf\|u_0^k\|_{L^3(R^3)}\leq\rho_{\max}, 
\end{equation}
by the definition of $\rho_{\max}$, we must have $\|u_0\|_{L^3(R^3)}=\rho_{\max}$. Thus we have $u_0^k\rightharpoonup u_0$ in $L^3$ and $\|u_0^k\|_{L^3}\to \|u_0\|_{L^3}$. Since $L^3(R^3)$ is uniformly convex as a Banach space, this also implies $u_0^k$ converges strongly to $u_0$ in $L^3(R^3)$.
The theorem is proved.\\

\medskip

From results above the following corollary can be proved following the same arguments as in \cite{RuSv}:
\begin{corollary}
Assume that every solution of the Cauchy problem \ref{eq:no1} with $u_0\in L^3(R^3)$ is regular, i.e. $T_{\max}=+\infty$ for each $u_0\in L^3(R^3)$. Then for $l=0,1,2,\dots$ there exist functions $F_l:[0,\infty)\to [0,\infty)$ such that
\begin{equation}
 t^{(l+1)/2}\sup_x|\nabla^lu(x,t)|\leq F_l(\|u_0\|_{L^3})\quad\mbox{for all $t>0$.}
\end{equation}

\end{corollary}
\end{section}

\end{document}